\documentclass{article}

\usepackage{arxiv}

\usepackage[utf8]{inputenc} 
\usepackage[T1]{fontenc}    
\usepackage{hyperref}       
\usepackage{url}            
\usepackage{booktabs}       
\usepackage{amsfonts}       
\usepackage{nicefrac}       
\usepackage{microtype}      
\usepackage{lipsum}

\usepackage{url}

\usepackage{amsthm,amsfonts,amssymb,amsmath,amsbsy,amsthm}
\usepackage{graphicx}
\usepackage{color}
\usepackage{algorithm,algorithmic}
\usepackage{soul}
\usepackage{caption}
\usepackage{subcaption}
\usepackage{lipsum} 
\usepackage{todonotes}
\usepackage{comment}

\newcommand{\bea}{\begin{eqnarray}}
\newcommand{\eea}{\end{eqnarray}}

\newcommand{\be}{\begin{equation}}
\newcommand{\ee}{\end{equation}}

\newcommand{\ben}{\begin{equation*}}
\newcommand{\een}{\end{equation*}}

\newcommand{\bean}{\begin{eqnarray*}}
\newcommand{\eean}{\end{eqnarray*}}

\DeclareMathOperator*{\argmin}{arg\,min}

\title{Quantifying the Multi-Objective Cost of Uncertainty}

\author{
  Byung-Jun Yoon\\
  Department of Electrical and Computer Engineering\\
  Texas A\&M University\\
  College Station, TX 77843, USA\\
  \texttt{bjyoon@ece.tamu.edu} \\
   \And
  Xiaoning Qian\\
  Department of Electrical and Computer Engineering\\
  Texas A\&M University\\
  College Station, TX 77843, USA\\
  \texttt{xqian@ece.tamu.edu} \\
    \And
  Edward R. Dougherty\\
  Department of Electrical and Computer Engineering\\
  Texas A\&M University\\
  College Station, TX 77843, USA\\
  \texttt{edward@ece.tamu.edu} \\
}

\begin{document}
\maketitle

\begin{abstract}
Various real-world applications involve modeling complex systems with immense uncertainty and optimizing multiple objectives based on the uncertain model. Quantifying the impact of the model uncertainty on the given operational objectives is critical for designing optimal experiments that can most effectively reduce the uncertainty that affect the objectives pertinent to the application at hand. In this paper, we propose the concept of mean multi-objective cost of uncertainty (multi-objective MOCU) that can be used for objective-based quantification of uncertainty for complex uncertain systems considering multiple operational objectives. We provide several illustrative examples that demonstrate the concept and strengths of the proposed multi-objective MOCU. Furthermore, we present a real-world example based on the mammalian cell cycle network to demonstrate how the multi-objective MOCU can be used for quantifying the operational impact of model uncertainty when there are multiple, possibly competing, objectives.\footnote{This work has been submitted to the IEEE for possible publication. Copyright may be transferred without notice, after which this version may no longer be accessible.}
\end{abstract}



\section{Introduction}

Investigating real-world systems and phenomena typically requires complex models that involve a large number of parameters. Even with sizeable amount of observation data, the high complexity of the model may render accurate parameter estimation impossible. While finding a reliable point estimate of the parameter vector may not be possible in such a case, it may be possible to identify the parameter ranges based on the available data and/or prior system knowledge, or in a more general setting, we may assume a joint distribution of the model parameters. Since different parameter values are possible, this gives rise to an uncertainty class of all possible models~\cite{Yoon2013tsp, Dougherty2019objectiveUQ}. Furthermore, this naturally places the uncertain model in a Bayesian framework, where the likelihood of every possible model in the uncertainty class is described by a prior distribution that could be constructed from prior system knowledge and/or existing data. For example, the MKDIP (maximal knowledge-driven information priors) proposed in~\cite{Boluki2017mkdip_bmc, Boluki2017mkdip_tcbb} shows how relational knowledge between interacting variables in the model can be used to construct the prior through constrained optimization based on a general framework of constraints stated in the form of conditional probabilities. The MKDIP technique has been previously shown to effectively translate biological pathway knowledge into a prior distribution for Bayesian learning, especially, for Bayesian classification and regression~\cite{Boluki2017mkdip_bmc, Boluki2017mkdip_tcbb}.

Given an uncertain model and its uncertainty class, how can one mathematically quantify the amount of uncertainty present in the model? Common approaches include estimating the variance or entropy of the uncertain parameters, as they both provide a simple and intuitive measure of the model uncertainty. However, they both have a critical downside from a practical perspective. In practical applications that involve mathematical modeling of a complex system, one cares about the model as it can serve as a vehicle for designing an effective operator (\textit{i.e.}, controller, classifier, filter) that can act on the system of interest or the data produced therefrom. As a result, if the model uncertainty does not affect the performance of the operator--even if the variance or the entropy of the model parameters might be substantial--one may not be concerned about the uncertainty from a practical perspective~\cite{Yoon2013tsp}.

Objective uncertainty quantification (objective-UQ)~\cite{Dougherty2019objectiveUQ} effectively addresses this shortcoming by measuring the impact that the model uncertainty has on the operator performance by estimating the mean objective cost of uncertainty (MOCU), originally proposed in~\cite{Yoon2013tsp}. MOCU quantifies the expected increase of the operational cost due to the uncertainty of the model. The uncertainty present in the model results in a potential operational cost increase as it prevents one from designing the operator that is optimal for the underlying true model (which is unknown) and necessitates a robust operator that guarantees good overall performance for all potential models included in the uncertainty class. While using a robust operator can guarantee good expected performance no matter what the actual underlying model is, it is generally suboptimal for any specific model in the uncertainty class, hence resulting in a cost increase.

As MOCU enables objective-based uncertainty quantification, it  provides an effective means of quantifying the expected impact of potential experiments on reducing model uncertainty that directly affects operator performance. For this reason, MOCU has been recently utilized in various application domains for optimal experimental design (OED), where examples include controlling uncertain gene regulatory networks (GRNs)~\cite{Dehghannasiri15bmc,Dehghannasiri15tcbb}, synchronization of uncertain Kuramoto oscillator models~\cite{Hong2021}, designing materials with targeted functional properties~\cite{Talapatra2018}, optimal sequential sampling~\cite{Broumand2015pr}, and active learning for optimal Bayesian classification~\cite{Zhao2021SMOCU,Zhao2021WMOCU}.

Originally, MOCU~\cite{Yoon2013tsp} was proposed for applications with a single objective, while many real-world applications require taking multiple objectives into account. For example, in drug design, one aims to find drug candidates with the greatest therapeutic efficacy but also with the smallest potential of drug-induced liver injury (DILI)~\cite{Agarwal2010DILI}. In materials discovery, one may be interested in designing new compounds that enhance multiple target properties~\cite{Gopakumar2018}, where some properties may even compete against each other. While various techniques exist for multi-objective optimization, no existing technique can be used for objective-based quantification of model uncertainty to the best of our knowledge.

In this paper, we extend the definition of MOCU for uncertain complex systems in case of multiple operational goals. This work will pave the way for multi-objective OED for various real-world applications--including drug design and materials discovery--that involve multi-objective optimization based on models with substantial uncertainty.
The paper is organized as follows. In Sec.~\ref{sec:mocu}, we provide a brief review of the single-objective MOCU and extend it for the case when there are two objectives, accompanied by additional insights and motivation for the proposed extension. Section~\ref{sec:simulation} presents simulation results that demonstrate the efficacy of the multi-objective MOCU for capturing the operational impact of model uncertainty on multiple objectives. A real-world example based on the mammalian cell cycle network is provided in Sec.~\ref{sec:grn}, focusing on objective-based uncertainty quantification for robust structural intervention with multiple intervention goals. In Sec.~\ref{sec:discussion}, we conclude the paper with further discussions on the significance of the proposed multi-objective MOCU, relation to other existing methods, and  potential future applications.


\section{Objective Uncertainty Quantification}
\label{sec:mocu}

In this section, we first review the definition of single-objective MOCU, which we originally proposed in~\cite{Yoon2013tsp}, and then extend the definition to allow objective-based uncertainty quantification under multi-objective settings.

\subsection{Brief review of single-objective MOCU}
\label{sec:mocu_single}

Let $\theta\in\Theta$ be the parameter vector of the uncertain model that belongs to an uncertainty class $\Theta$ comprised of all possible models. We denote $\psi \in \Psi$ as an operator in an operator class $\Psi$, where $\xi_\theta(\psi)$ is the cost of applying the operator $\psi$ when the true model is $\theta$. Given full knowledge of the model $\theta$, we can design the optimal operator $\psi_\theta$ that minimizes the operational cost as follows
\be
    \psi_\theta = \argmin_{\psi\in\Psi} \xi_\theta(\psi).    \label{eq:opt_operator}
\ee
When the true model $\theta$ is not precisely known, the optimal operator for the true underlying model cannot be designed, in which case it is desirable to design a robust operator that keeps the expected operational cost at a minimum for all possible models $\theta\in\Theta$. For example, we can design the optimal robust operator
\be
    \psi^* = \argmin_{\psi\in\Psi} E_\Theta \Big [\xi_\theta(\psi) \Big ]
\ee
that minimizes the expected operational cost $\xi_\theta(\psi)$, where the expectation is taken with respect to the prior distribution $\pi(\theta)$ of the model $\theta\in\Theta$. The prior $\pi(\theta)$ may be mathematically constructed from the available prior knowledge regarding the system being modeled, for which techniques such as the MKDIP (Maximal Knowledge-Driven Information Priors)~\cite{Boluki2017mkdip_bmc,Boluki2017mkdip_tcbb} can be used. If no such prior knowledge is available, one may simply use a noninformative uniform prior. 
MOCU~\cite{Yoon2013tsp} quantifies the objective cost of uncertainty as follows:
\be
    \eta = E_\Theta \Big [\xi_\theta(\psi^*) - \xi_\theta(\psi_\theta) \Big ]  \label{eq:mocu-so}
\ee
by measuring the increase of operational cost due to using a robust operator rather than the optimal operator for the true model, which is inevitable since the true model is unknown.


\subsection{Objective uncertainty quantification for two objectives}
\label{sec:mocu_two}

Now suppose that we have two different cost functions of interest, $\xi^1_\theta(\psi)$ and $\xi^2_\theta(\psi)$, which measure the operational performance of a given operator $\psi$ for two different objectives, respectively. How can we quantify the objective cost of uncertainty in this case? A reasonable approach would be to define a weighted cost function as follows
\be
    \xi_\theta(\psi, \lambda) = \lambda  \xi^1_\theta(\psi) + (1-\lambda) \xi^2_\theta(\psi)  \label{eq:weighted_cost}
\ee
where $\lambda \in [0,1]$ is a weight parameter. Based on this weighted cost function, one could assess to what extent the model uncertainty increases the combined cost for the two objectives at hand. In fact, similar approaches are often used in multi-objective optimization problems, as it reduces multiple objective functions into a single weighted objective function, thereby keeping the optimization problem computationally tractable.

For any given $\lambda \in [0,1]$, we can define the optimal model-specific operator as follows
\be
    \psi_\theta(\lambda) = \argmin_{\psi\in\Psi} \xi_\theta(\psi, \lambda), \label{eq:psi_theta_lambda}
\ee
which minimizes the combined cost function $\xi_\theta(\psi, \lambda)$ for the model $\theta$. Similarly, we can also define the optimal robust operator
\be
    \psi^*(\lambda) = \argmin_{\psi\in\Psi} E_\Theta \Big [\xi_\theta(\psi, \lambda) \Big ],  \label{eq:psi*_lambda}
\ee
which minimizes the expected value of the combined cost function for all possible models $\theta \in \Theta$. 
Based on (\ref{eq:psi_theta_lambda}) and (\ref{eq:psi*_lambda}), we can define the MOCU for a specific value of $\lambda$ as a function of $\lambda$ as follows:
\be
    \eta(\lambda) = E_\Theta \Big [\xi_\theta\big(\psi^*(\lambda), \lambda\big) - \xi_\theta\big(\psi_\theta(\lambda), \lambda\big) \Big ]. \label{eq:mocu_lambda}
\ee
It should be noted that $\eta(\lambda)\geq 0$ for any $\lambda \in [0,1]$. When $\lambda=1$, $\eta(\lambda)$ becomes the single-objective MOCU for using the first cost function $\xi^1_\theta(\psi)$. When $\lambda=0$, $\eta(\lambda)$ becomes identical to the second single-objective MOCU for using the cost function $\xi^2_\theta(\psi)$.

It is critical to remember that $\eta(\lambda)$ given by \eqref{eq:mocu_lambda} cannot be simply obtained by computing the weighted average of $\eta(1)$ (i.e., the single-objective MOCU based on the first objective) and $\eta(0)$ (i.e., the single-objective MOCU based on the second objective). The optimal robust operator $\psi^*(\lambda)$ depends on the value of $\lambda$ in a highly complex manner, and by no means can it be obtained by interpolating the robust operators for different values of $\lambda$ -- \textit{e.g.}, via interpolation of $\psi^*(0)$ and $\psi^*(1)$.
Consequently, the $\lambda$-specific MOCU $\eta(\lambda)$ in~\eqref{eq:mocu_lambda} cannot be simply obtained from a linear combination of single-objective MOCU values (\textit{i.e.}, $\eta(1)$ and $\eta(0)$), as $\eta(\lambda)$ depends on the optimal robust operator $\psi^*(\lambda)$ as well as the optimal model-specific operator $\psi_\theta(\lambda)$.

Now, a natural question arises: how should we select the value of $\lambda$?  When a weighted cost function -- like the one in \eqref{eq:weighted_cost} -- is used for multi-objective optimization, the choice is often somewhat arbitrary. Choosing a specific value for $\lambda$ requires predetermining the relative importance of the two cost functions. Instead of choosing an arbitrary $\lambda$, we can compute the average cost of uncertainty for all possible values of $\lambda$ in the following way
\be
    \eta = \int_0^1 \eta(\lambda) d\lambda.     \label{eq:mocu_two_ave}
\ee
We can generalize this further by estimating the expectation of $\eta(\lambda)$ with respect to a distribution $p(\lambda)$ of $\lambda$ as follows
\be
    \eta = E_{\lambda} \Big[ \eta(\lambda) \Big] = \int \eta(\lambda)p(\lambda) d\lambda,     \label{eq:mocu_two_mean}
\ee
where \eqref{eq:mocu_two_ave} is a special case when $p(\lambda)$ is a uniform distribution for $\lambda \in [0,1]$.
This double-objective MOCU
extends the single-objective MOCU given by~\eqref{eq:mocu-so} in a nontrivial manner by adding novel dimensions to the original definition: (i) the optimization of the robust operator for two (possibly competing) objectives and (ii) the quantification of the expected cost of model uncertainty by assessing the performance degradation of this robust operator that is optimized for both objectives simultaneously.


\subsection{Multi-objective MOCU}
\label{sec:mocu_multi}

We can further generalize the definition of MOCU for two objectives, shown in \eqref{eq:mocu_two_mean}, for the cases when we have multiple objectives. Suppose there are $n$ different cost functions 
\be
    \xi^1_\theta(\psi), \, \xi^2_\theta(\psi), \, \cdots, \, \xi^n_\theta(\psi),
\ee
where $\xi^k_\theta(\psi)$ is the cost function that is used to assess the operational cost of a given operator $\psi$ with respect to the $k$-th objective at hand. Let $\boldsymbol{\lambda} = (\lambda_1, \lambda_2, \cdots, \lambda_n)$ be a weight vector such that $\sum_{i=1}^{n} \lambda_i = 1$ and $\lambda_i \geq 0$. Given $\boldsymbol{\lambda}$, we define the following combined cost function as follows
\be
    \xi_\theta(\psi, \boldsymbol{\lambda}) = \sum_{i=1}^{n} \lambda_i  \xi^i_\theta(\psi),
\ee
which gives us the weighted operational cost of applying the operator $\psi$ to the model $\theta$.
For any given $\boldsymbol{\lambda}$, we can define the optimal model-specific operator as follows:
\be
    \psi_\theta(\boldsymbol{\lambda}) = \argmin_{\psi\in\Psi} \xi_\theta(\psi, \boldsymbol{\lambda}).  \label{eq:psi_theta_lambda_multi}
\ee
Similarly, we can also define the optimal robust operator:
\be
    \psi^*(\boldsymbol{\lambda}) = \argmin_{\psi\in\Psi} E_\Theta \Big [\xi_\theta(\psi, \boldsymbol{\lambda}) \Big ].  \label{eq:psi*_lambda_multi}
\ee
As in Sec.~\ref{sec:mocu_two}, based on (\ref{eq:psi_theta_lambda_multi}) and (\ref{eq:psi*_lambda_multi}), we can define the MOCU for a specific value of $\boldsymbol{\lambda}$ as a function of $\boldsymbol{\lambda}$ as follows:
\be
    \eta(\boldsymbol{\lambda}) = E_\Theta \Big [\xi_\theta(\psi^*(\boldsymbol{\lambda}), \boldsymbol{\lambda}) - \xi_\theta(\psi_\theta, \boldsymbol{\lambda}) \Big ]. \label{eq:mocu_lambda_multi}
\ee
As before, in order to avoid choosing a fixed weight vector $\boldsymbol{\lambda}$, we instead define the multi-objective MOCU (i.e., the mean multi-objective cost of uncertainty) as follows:
\be
    \eta_{\scriptsize{\mbox{multi}}} =  E_{\boldsymbol{\lambda}} \Big[ \eta(\boldsymbol{\lambda}) \Big] =
    \int \eta(\boldsymbol{\lambda}) p(\boldsymbol{\lambda}) \ d\boldsymbol{\lambda}, \label{eq:mocu_multi}
\ee
where $p(\boldsymbol{\lambda})$ is the distribution of the weight vector. When we do not have a specific preference regarding $\boldsymbol{\lambda}$, we can use a uniform distribution for $p(\boldsymbol{\lambda})$ such that $\boldsymbol{\lambda}$ is uniformly distributed on the hyperplane that satisfies $\lambda_i \geq 0$ for $\forall i$ and $\sum_{i=1}^{n} \lambda_i = 1$, which is simply a flat Dirichlet distribution.


\section{Simulation Results}
\label{sec:simulation}

In this section, we consider a multi-objective optimization problem under uncertainty to demonstrate the efficacy of the multi-objective MOCU proposed in Sec.~\ref{sec:mocu}. For this purpose, we define two objective functions $f_1(x,y)$ and $f_2(x,y)$, where the first objective function is defined as
\be
    f_1(x,y) = \alpha_1 (x-\gamma_1)^2 + \beta_1 (y-\delta_1)^2 \label{eq:obj1}
\ee
and the second objective is defined as
\be
    f_2(x,y) = \alpha_2 (x-\gamma_2)^2 + \beta_2 (y-\delta_2)^2.    \label{eq:obj2}
\ee
The goal of this multi-objective optimization is to minimize both objectives, where some (or all) of the parameters ({\it i.e.}, $\alpha_i$, $\beta_i$, $\gamma_i$, $\delta_i)$) may be uncertain. The functions \eqref{eq:obj1} and \eqref{eq:obj2} generalize test functions proposed in~\cite{Binh1997, Schaffer1985}, which are frequently used for testing multi-objective optimization algorithms.

\begin{figure}[t!]
    \centering
    \includegraphics[width=4in]{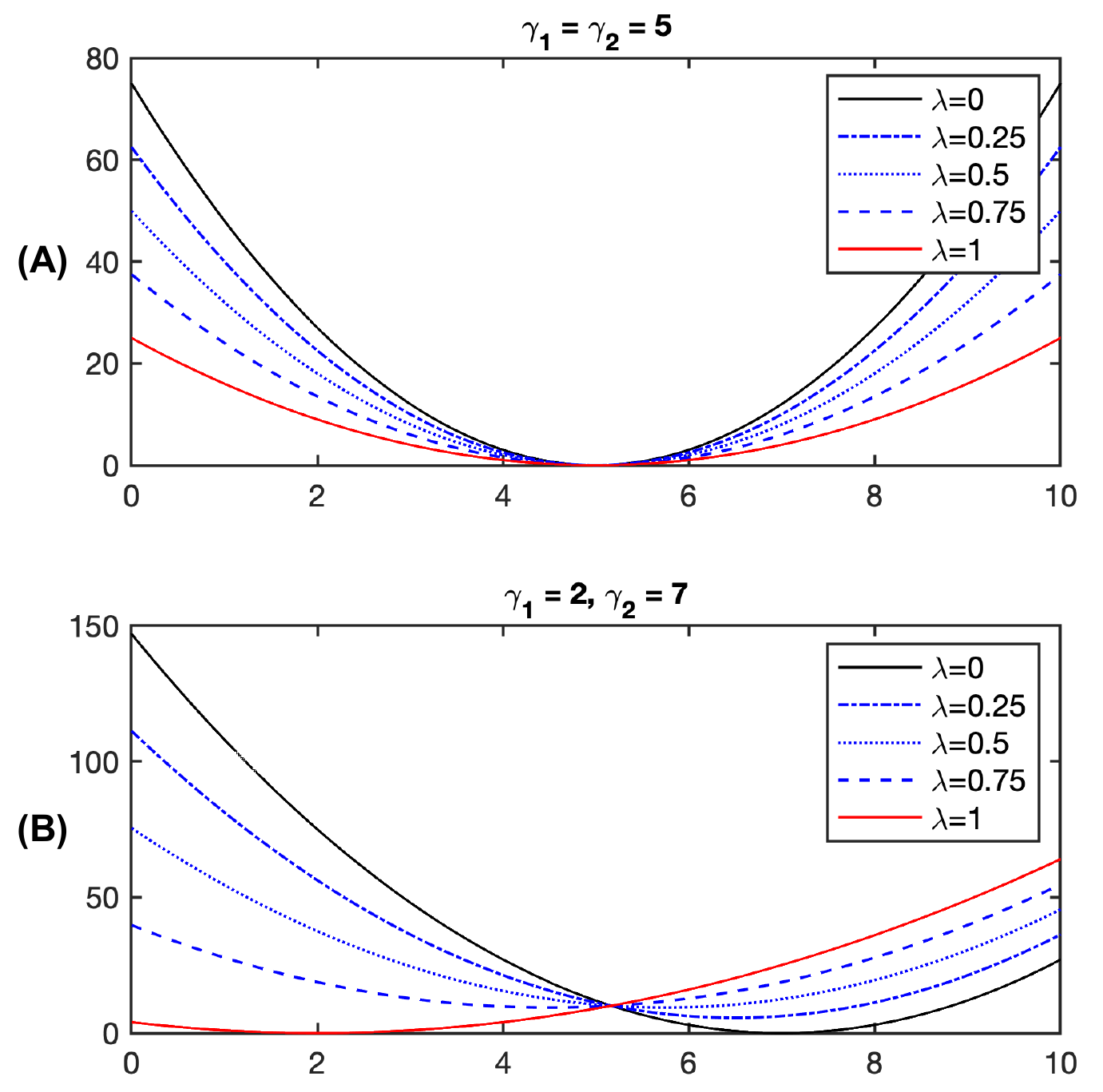}
\caption{Illustrative examples. (A) The global minimum does not depend on $\lambda$ nor $\alpha_1, \alpha_2$, hence the model uncertainty does not affect achieving the objective. (B) The location of the global minimum changes with $\lambda, \alpha_1, \alpha_2$, as a result of which, the uncertainty induces an additional operational cost.}
    \label{fig:example}
\end{figure}

How does the presence of uncertainty in the parameters affect the objectives--\textit{i.e.}, minimization of the $f_1(x,y)$ and $f_2(x,y)$? To answer this question, we consider two simple examples. Let $g_1(x) = \alpha_1 (x-\gamma_1)^2$ and $g_2(x) = \alpha_2 (x-\gamma_2)^2$, where we aim to minimize  $g(x,\lambda) = \lambda g_1(x) + (1-\lambda) g_2(x)$. Figure~\ref{fig:example}A shows the objective function $g(x,\lambda)$ for the case when $\gamma_1 = \gamma_2 = 5$, $\alpha_1 = 1$, $\alpha_2 = 3$ for various values of $\lambda$. As we can see in Fig.~\ref{fig:example}A, as $\gamma_1$ and $\gamma_2$ are identical, different values of $\lambda$ does not affect the location of the global minimum. This is true for any $\alpha_1 \geq 0$ and $\alpha_2 \geq 0$, which clearly shows that uncertainty regarding $\alpha_1$ and $\alpha_2$ does not affect the operational goal of minimizing $g(x,\lambda)$. 
Now, let us consider the case when $\gamma_1 = 2$ and $\gamma_2 = 7$, while $\alpha_1$ and $\alpha_2$ remain the same. Figure~\ref{fig:example}B depicts $g(x,\lambda)$ in this case for different values of $\lambda$. We can see that the global minimum changes for different $\lambda$. In fact, different values of $\alpha_1$, $\alpha_2$, and $\lambda$ will affect the location of the global minimum in this case, which immediately shows why any uncertainty in $\alpha_1$ and $\alpha_2$ would affect the operational goal in this case. 

While the proposed multi-objective MOCU can effectively capture the impact of such model uncertainty on the objective and quantify the ``objective'' cost of uncertainty, traditional measures such as entropy and variance fail to do so as they are not designed to measure the uncertainty pertinent to a specific objective.

\subsection{Case-1}
\label{sec:case1}

\begin{figure}[t!]
    \centering
    \includegraphics[width=4in]{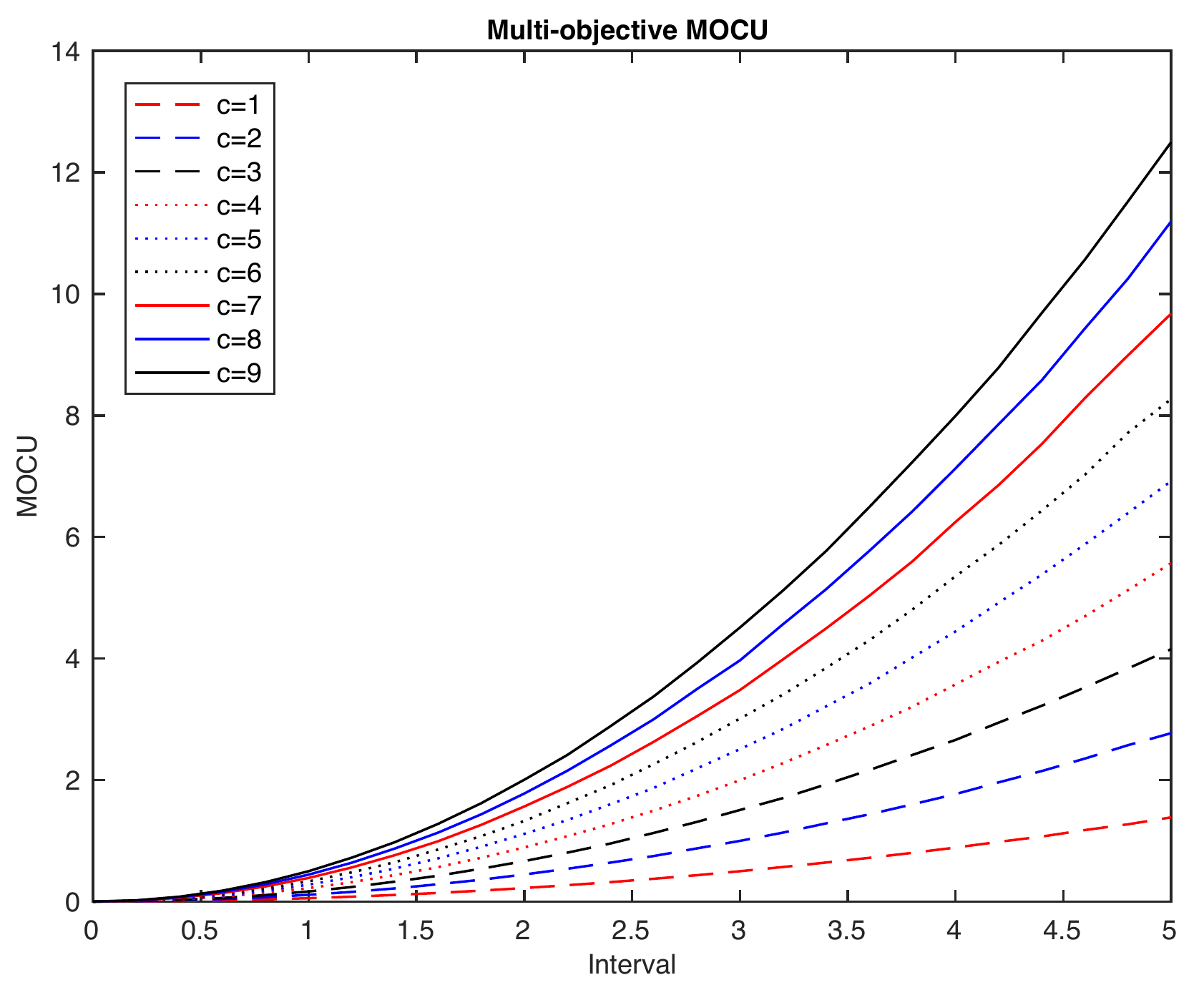}
\caption{Multi-objective MOCU estimation for Case-1. Increase of the interval ($\Delta$) of the uncertain parameters increases the objective uncertainty measured by MOCU. For a given interval, the value of $c$ affects how the model uncertainty impacts the operational cost.}
    \label{fig:mocu1}
\end{figure}

Here, we consider the two objective functions in \eqref{eq:obj1} and \eqref{eq:obj2}, where $\alpha_1 = \alpha_2 = \beta_1 = \beta_2 = c$, $\gamma_1 = \delta_1 = 0$, and $\gamma_2, \delta_2 \in [0,\Delta]$. Uncertainty is present in the parameters $\gamma_2$ and $\delta_2$, both of which are uniformly distributed in $[0,\Delta]$, while all other parameters are known. Figure~\ref{fig:mocu1} shows the multi-objective MOCU estimated by \eqref{eq:mocu_two_ave} for various values of $c$. As shown in Fig.~\ref{fig:mocu1}, objective uncertainty (measured by MOCU) increases as the uncertain interval ($\Delta$) increases. The graphs also show that the uncertainty impacts the performance more significantly when $c$ is larger. In comparison, while both entropy and variance will increase with $\Delta$, neither of them can capture the impact of varying $c$, as they are constant given the uncertain interval, independent of different values of $\alpha_1$, $\alpha_2$, $\beta_1$, $\beta_2$.

\subsection{Case-2}
\label{sec:case2}

\begin{figure}[t!]
    \centering
    \includegraphics[width=4in]{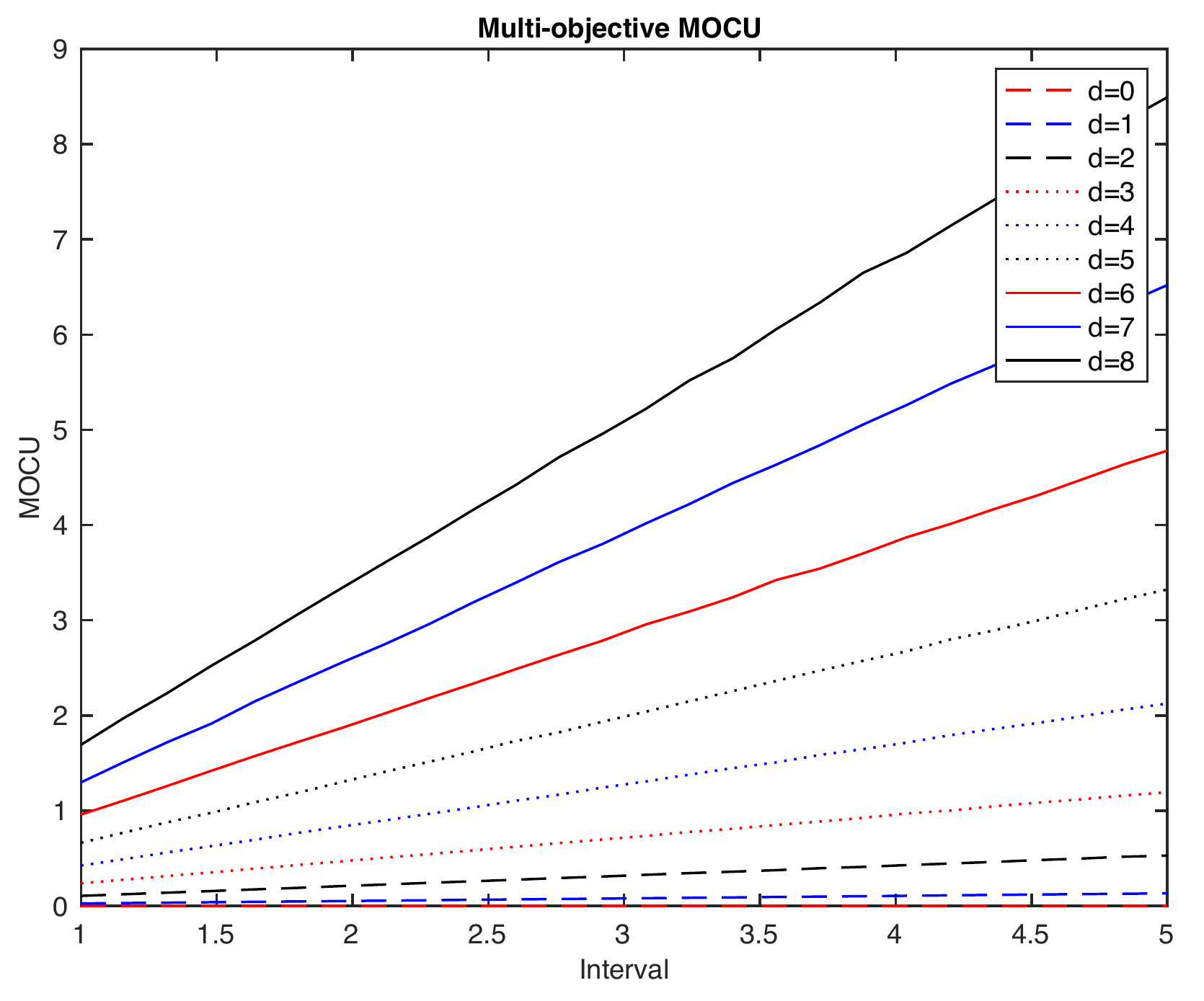}
\caption{Multi-objective MOCU estimation for Case-2. Increasing the interval ($\Delta$) leads to an increase of the objective uncertainty measured by MOCU. The value of $d$ affects how the model uncertainty impacts the operational cost. For $d=0$, we can see that the model uncertainty does not increase the cost at all.}
    \label{fig:mocu2}
\end{figure}

Suppose that the parameters $\alpha_1, \alpha_2, \beta_1, \beta_2$ are uncertain and that they are uniformly distributed with $\alpha_1, \alpha_2, \beta_1, \beta_2 \in  [0,\Delta]$. The other parameters are known, where $\gamma_1 = \delta_1 = 0$ and $\gamma_2 = \delta_2 = d$. Figure~\ref{fig:mocu2} shows the estimated multi-objective MOCU as a function of $\Delta$ ({\it i.e.}, length of the uncertain interval) for the unknown parameters. As expected, MOCU increases with a larger $\Delta$. The impact of this uncertainty is larger when $\gamma_2 - \gamma_1$ and $\delta_2 - \delta_1$ tends to be larger ({\it i.e.}, for larger $d$), which is intuitive when we examine \eqref{eq:obj1} and \eqref{eq:obj2}. One interesting point to note is that the MOCU is zero when $d=0$, implying that the uncertainty present in the parameters $\alpha_1, \alpha_2, \beta_1, \beta_2$ does not impact the objective at all. It should be noted that neither the entropy nor the variance can be used to predict that the uncertainty present in \eqref{eq:obj1} and \eqref{eq:obj2} will not impact the objective when $d=0$. Furthermore, the entropy and the variance are unable to quantify the impact of $d$ on the objective, since they will only be affected by $\Delta$ and not by $d$.

\subsection{Case-3}
\label{sec:case3}

\begin{figure}[t!]
    \centering
    \includegraphics[width=4in]{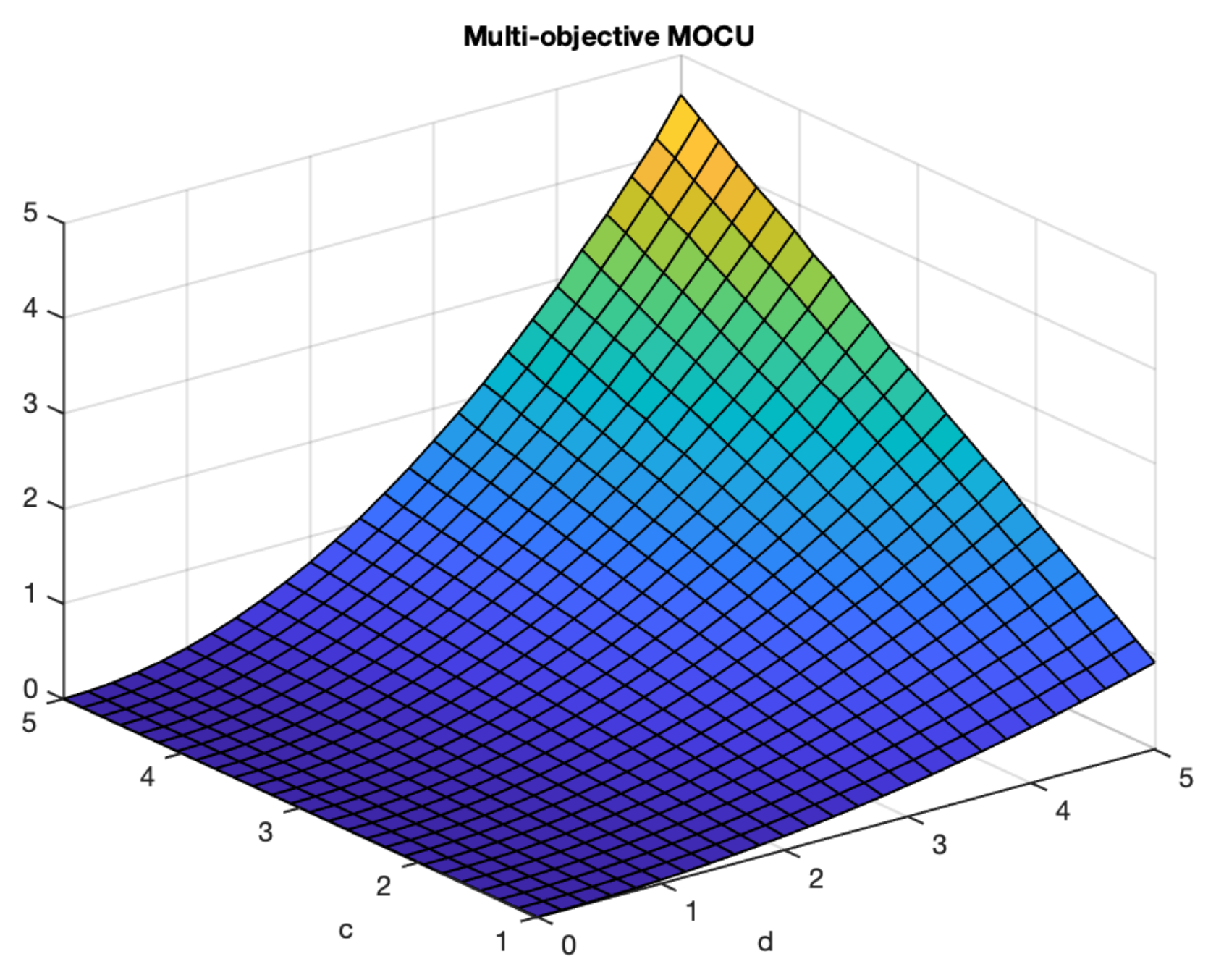}
\caption{Multi-objective MOCU estimation for Case-3.}
    \label{fig:mocu3}
\end{figure}

Finally, we consider the case when $\alpha_1, \alpha_2, \beta_1, \beta_2$ are uncertain and uniformly distributed in $[0,c]$, $\gamma_1 = \delta_1 = 0$, and $\gamma_2, \delta_2$ are uncertain and uniformly distributed in $[0,d]$. From the results in Sec.~\ref{sec:case1}, we already know how the uncertainty regarding $\gamma_2, \delta_2 \in [0,d]$ affects the objectives. Similarly, results in Sec.~\ref{sec:case2} show how the uncertainty $\alpha_1, \alpha_2, \beta_1, \beta_2 \in [0,c]$ would impact the objectives. While measures such as the entropy and variance are unable to quantify the impact of uncertainty on the objectives, the proposed multi-objective MOCU can effectively quantify the operational cost increase resulting from the uncertainty present in the objective functions. Figure~\ref{fig:mocu3} shows the estimated MOCU as a function of $c$ and $d$ for $1\leq c \leq 5$ and $0 \leq d \leq 5$.

\section{Quantifying the Multi-Objective Cost of Intervention in Uncertain Gene Regulatory Networks}
\label{sec:grn}

\begin{figure}[b!]
    \centering
    \includegraphics[width=4in]{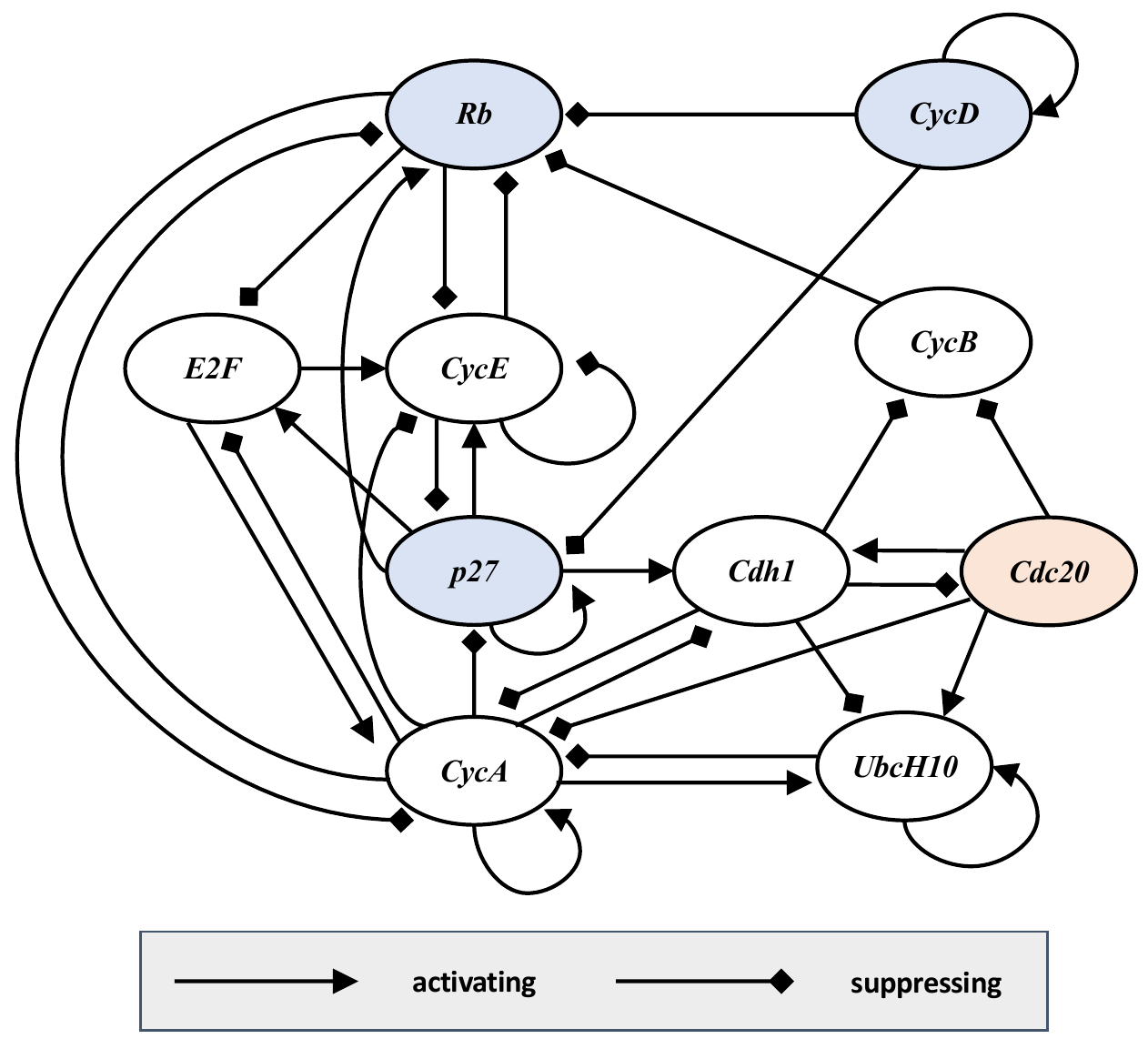}
\caption{Illustration of the mammalian cell cycle network. The operational goal is to desirably change the network dynamics via robust structural intervention in the presence of model uncertainty. The first objective is to minimize the steady-state probability mass in undesirable states (that belong to the set $U$), where the genes \textit{CycD}, \textit{Rb}, \textit{p27} (shown in light blue) are simultaneously down-regulated. The second objective is to minimize the steady-state probability mass of potentially pathological states (that belong to the set $P$), in which \textit{Cdc20} (shown in light red) is up-regulated.}
    \label{fig:mamccl}
\end{figure}

In this section, we provide a real-world example based on the Mammalian Cell Cycle pathway~\cite{Faure2006} that demonstrates how the proposed multi-objective MOCU can be utilized for objective-based uncertainty quantification of complex systems when there are multiple operational objectives.
The gene regulatory network model for the Mammalian Cell Cycle pathway~\cite{Faure2006} consists of ten genes -- \textit{CycD}, \textit{Rb}, \textit{p27}, \textit{E2F}, \textit{CycE}, \textit{CycA}, \textit{Cdc20}, \textit{Cdh1}, \textit{UbcH10}, and \textit{CycB} -- whose regulatory relationships (i.e., activating or suppressing) are illustrated in Figure~\ref{fig:mamccl}. As in~\cite{Yoon2013tsp}, a Boolean Network with perturbation~(BNp)~\cite{PBNBook} is adopted to model the gene expression dynamics, denoted by $\mathbf{x}(t) = [x_1(t), x_2(t), \cdots, x_{10}(t)]$ in the aforementioned order of ten genes. Based on the majority voting rule with the regulatory relationships between these genes, the gene expression will be turned on ($x_i=1$) or off ($0$) depending on whether the majority of the expressed regulators are activators or suppressors. In this BNp model, we set the random perturbation probability to $p=0.01$ so that each gene has a small probability of $p$ to have a randomly flipped gene expression state at each time point.

Mammalian cell division is controlled via extra-cellular signals and it is coordinated with overall cellular growth. The positive signals, or growth factors, instigate \textit{CycD}. The cyclins inhibit \textit{Rb} by phosphorylation. The gene \textit{p27} can stop the uncontrolled cell cycle as it blocks the action of \textit{CycE} or \textit{CycA} and thereafter \textit{Rb} can also be expressed, even in the presence of \textit{CycE} or \textit{CycA}. When \textit{CycD}, \textit{Rb}, and \textit{p27} are simultaneously down-regulated (i.e., $x_1=x_2=x_3=0$), the cell cycles can continue indefinitely even in the absence of any growth factor, representing cancerous phenotypes. Therefore, we define such network states as \textit{undesirable} states $U=\{\mathbf{x}|x_1=x_2=x_3=0\}$. The aim of therapeutic intervention is then to reduce the probability of entering into these undesirable states in $U$. In particular, we focus on the class of structural intervention strategies~\cite{Qian2008,Qian2011} that block the regulatory action between any pair
of genes in the given network in Figure~\ref{fig:mamccl}. Here in our experiments, we are constrained to block only one regulation while more flexible intervention can be derived in a straightforward manner as detailed in~\cite{Qian2008,Qian2011}. The intervention objective is to minimize the steady-state mass of these undesirable states:
\be
    \min \{\pi_U = \sum_{\forall \mathbf{x} \in U} \pi_{\mathbf{x}} \},
\ee
where $\pi_{\mathbf{x}}$ denotes the steady-state probability of the corresponding network state $\mathbf{x}$. This is our primary objective for deriving the optimal robust intervention, typically referred to as the IBR (intrinsically Bayesian robust) structural intervention~\cite{Yoon2013tsp}.

As discussed in~\cite{Qian2012}, in addition to reducing the undesirable steady-state mass $\pi_U$, it is also critical to make sure that the intervention applied to the network does not incur unforeseen collateral damage. Hence, when designing intervention strategies for  gene regulatory network models, we may want to further constrain the shifted steady-state mass to network states that are known to be ``safe'' and do not correspond to pathological phenotypes. In order to achieve this goal, we can further incorporate prior knowledge and strive to minimize the steady-state mass of phenotypically constrained states $P$. In this set of experiments, we define $P = \{\mathbf{x}|x_7=1\} \setminus U$, where $\setminus$ denote the set difference operator. In fact, $x_7=1$ denotes the network states with expressed \textit{Cdc20}, whose overexpression has been reported in breast, lung, gastric, and pancreatic cancers~\cite{Curtis2020}. Therefore, our second objective for deriving the IBR structural intervention is
\be
    \min\{\pi_P = \sum_{\forall \mathbf{x} \in P} \pi_{\mathbf{x}}\}.
\ee

To illustrate the effectiveness of objective-based uncertainty quantification using  multi-objective MOCU proposed in this paper, we assume that the regulatory information (either activating or suppressing) is unknown for some of the edges in Figure~\ref{fig:mamccl}, which gives rise to an uncertainty class of possible network models. We vary the amount of uncertainty in the mammalian cell cycle network model by varying the number of unknown edges from 1 to 8 and randomly sample the corresponding number of edges to form the uncertainty class $\mathbf{\Theta}$. This uncertainty class $\mathbf{\Theta}$ contains the network models with all possible combinations of regulatory relationships for the set of corresponding edges (i.e., edges with unknown regulatory relationships). We compute the multi-objective MOCU based on \eqref{eq:mocu_two_mean} for this uncertainty class to investigate how the quantified multi-objective uncertainty changes as a function of the number of edges in the network with unknown regulatory relationships.

\begin{figure}[t!]
    \centering
    \includegraphics[width=4in]{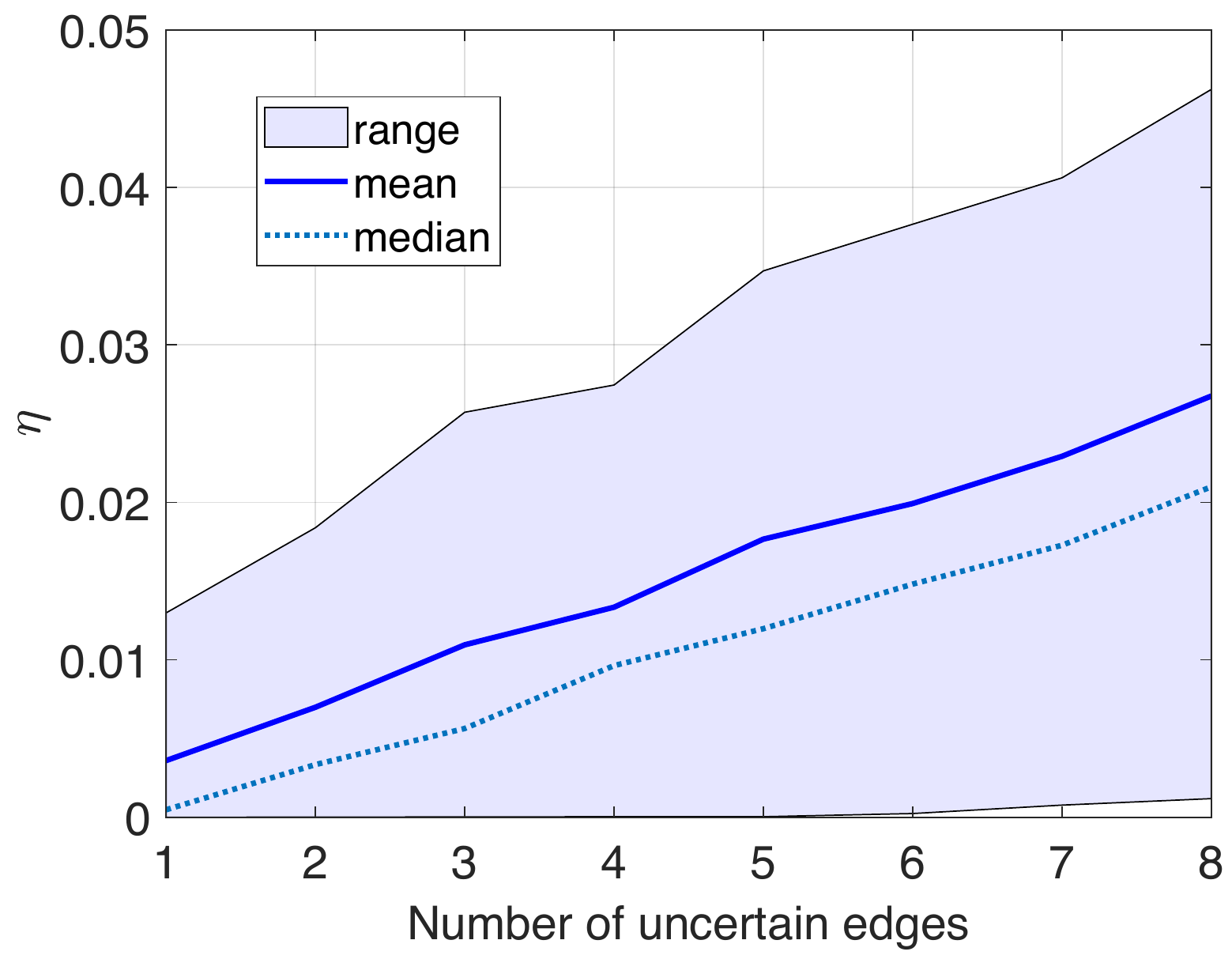}
\caption{Average multi-objective MOCU is shown as a function of the number of unknown regulatory relationships based on 500 random runs in the mammalian cell cycle network. The shaded region shows the range of corresponding multi-objective MOCU values from the minimum to the mean+1 std (standard deviation) for these 500 runs. We note that the size of the uncertainty class increases exponentially with the number of uncertain edges. For example, the size of the uncertainty class with 7 uncertain edges is $\binom{35}{7}$ as there are 35 edges in the mammalian cell cycle BNp. Therefore it may need a significantly higher number of runs to have more smooth curves without significant fluctuation.}
    \label{fig:mammocu}
\end{figure}

Figure~\ref{fig:mammocu} shows the minimum, median, mean, and mean+1 std (standard deviation values) of multi-objective MOCU in~\eqref{eq:mocu_two_mean} from 500 randomly sampled uncertainty classes (with replacement) based on the mammalian cell cyle network represented by BNp. It is clear that our multi-objective MOCU faithfully captures the increasing uncertainty affecting the performance of the derived structural intervention as the number of uncertain edges increases. It also confirms that it is critical to capture the ``objective-based'' uncertainty using the proposed multi-objective MOCU when the intervention objectives are of the ultimate importance, as clearly shown in the figure, the range of multi-objective MOCU also increases with the number of unknown edges. In fact, the median multi-objective MOCU values are consistently below the mean values. This implies that there are certain edges whose regulatory information is more critical, and as a result, the missing information regarding such edges may degrade the intervention performance more significantly than the others. Hence, with limited resource and time, it is more desirable and also more promising to design experiments to understand these regulatory actions for the purpose of deriving effective therapeutic intervention strategies.


\section{Discussions}
\label{sec:discussion}

In this paper, we introduced the definition of multi-objective MOCU, which extends the original single-objective MOCU proposed in~\cite{Yoon2013tsp} and enables objective-based quantification of model uncertainty in applications with multiple objectives. As entropy and variance play critical roles in various learning and inference problems, MOCU provides the foundation for optimal design of experiments and robust operators for complex uncertain systems. The proposed multi-objective MOCU enables optimal experimental design for various real-world applications--such as drug design and material discovery--that involve modeling complex systems with substantial uncertainty and multi-objective optimization based on the uncertain model.

The need for robust optimization of controllers, classifiers, estimators, or  various other types of operators in the presence of uncertainty arises frequently in diverse science and engineering problems that deal with complex real-world systems. While the uncertainty of a model representing the complex system of interest may be described by the probability distribution of its parameter vector in a Bayesian framework, the potential impact of this model uncertainty may significantly differ for different operations. For example, if we consider the gene regulatory network (GRN) model for the mammalian cell cycle network that was presented in Sec.~\ref{sec:grn}, the absence of regulatory knowledge regarding a specific pair of genes in the network shown in Figure~\ref{fig:mamccl} may considerably degrade the structural intervention performance of minimizing the steady-state mass in $U$ while affecting the intervention performance of minimizing the stead-state mass in $P$ only marginally.
As a result, when quantifying the uncertainty of the model, it is critically important to assess the operational impact of the uncertainty in a way that takes all operational objectives at hand into account.

Traditional uncertainty measures, such as entropy and variance, cannot serve this purpose as they are operation agnostic, hence do not inform us of the impact of the model uncertainty on the operations to be performed and the objectives to be achieved thereby. While the concept of MOCU, originally proposed in~\cite{Yoon2013tsp}, provides an effective means of such \textit{objective-based UQ}, its original definition in \eqref{eq:mocu-so} assumed a single objective and could not be used for quantifying uncertainty when there are multiple -- and potentially competing -- operational goals.
The multi-objective MOCU proposed in this paper effectively addresses these limitations.

Considering that the multi-objective MOCU provides an effective means of  multi-objective UQ, it can be directly used to extend the existing MOCU-based optimal experimental design (OED) strategies~\cite{Dehghannasiri15bmc,Dehghannasiri15tcbb,Talapatra2018} to enable the prediction of optimal experiment that is expected to enhance the performance of multi-objective robust operators. For example, the OED strategies in~\cite{Dehghannasiri15bmc,Dehghannasiri15tcbb} aimed at improving the structural intervention performance in an uncertain gene regulatory network model to minimize the steady-state probability mass in undesirable states corresponding to a single pathological phenotype. As demonstrated by the example discussed in Sec.~\ref{sec:grn}, the proposed multi-objective MOCU can extend the OED strategies for multiple objectives, which may involve minimizing the likelihood of multiple pathological phenotypes and/or maximizing the shift of the steady-state mass towards more desirable states or states that are not associated with any known aberrant cell behavior. 

Last but not least, we would like to emphasize that both the single-objective MOCU~\cite{Yoon2013tsp} as well as the multi-objective MOCU proposed in this work are consider the ``model uncertainty'' with respect to the underlying system. As such, when MOCU is used for OED~\cite{Hong2021,Dehghannasiri15bmc,Dehghannasiri15tcbb} or active learning~\cite{Zhao2021SMOCU,Zhao2021WMOCU}, the acquisition of experimental results or data will be guided by their potential impact on reduction of this model uncertainty. This is fundamentally different from existing Bayesian optimization~(BO)~\cite{Jones1998EGO,Couckuyt2014,Guang2018arxiv} and knowledge gradient~(KG) strategies~\cite{Frazier2009,Yahyaa2015} (either single- or multi-objective problems), as BO and KG probabilistically model the operational objectives using surrogate models of the ``evaluation  metrics'' -- rather than modeling the underlying system and the uncertainty therein -- for example using Gaussian processes (GPs) that are commonly adopted.
In this work, we focus on how one can quantify the model uncertainty directly impacting multiple operational objectives, which cannot be achieved by either BO or KG.
One point that is worth mentioning is that it has been established in~\cite{Shahin2019} that, under certain conditions, MOCU-based OED becomes equivalent to KG.
In the context of multi-objective BO, the optimization goal is to approach an estimated Pareto front, for which the operational objectives are often approximated by GPs.
The proposed multi-objective MOCU can help measure how model uncertainty may influence the estimation of the Pareto front and thereafter the results of multi-objective BO. One interesting research question here is that how the OED performance will change with the number of available training samples for different modeling and OED strategies, which we leave for our future research.

The problem of optimization and robust decision making under uncertainty arises frequently in various science and engineering domains~\cite{Kochenderfer2015,Diwekar2020}, as many real-world complex systems cannot be accurately modeled or completely identified in practice~\cite{Yoon2013tsp}. Relevant problems are abundant across diverse disciplines, including the robust intervention in gene regulatory networks ~\cite{Dehghannasiri15bmc,Dehghannasiri15tcbb}, optimization of structural materials~\cite{Beck2012}, robust design and operation of chemical engineering processes~\cite{Wendt2002}, optimization of the economic and life‐cycle environmental performance of industrial processes~\cite{Sabio2014}, optimal target selection for metabolic engineering~\cite{Wang2006a,Wang2006b}, optimal wake steering strategies for reducing power losses due to aerodynamic interactions between turbines~\cite{Quick2017}, and robust power dispatching in modern power grid systems~\cite{Chen2019a,Chen2019b}, just to give a few representative examples. The multi-objective MOCU proposed in this paper can provide effective means of objective-based quantification of uncertainties in the aforementioned (and other similar) systems. Furthermore, it can enable optimal experimental design~\cite{Dehghannasiri15bmc,Dehghannasiri15tcbb,Hong2021} and active learning~\cite{Zhao2021SMOCU,Zhao2021WMOCU} for objective-driven uncertainty reduction given multiple operational goals.

\pagebreak


\bibliographystyle{ieeetr}
\bibliography{multi-objective-mocu-references.bib}

\end{document}